\documentclass[12pt]{article}

\newcommand{\version}{November 26, 2021}

\usepackage{graphicx}
\usepackage{pictex}
\usepackage{mathrsfs}
\usepackage{amssymb}
\usepackage{color}

\def\hop#1{\relax} 

\title{Proof of a conjecture on hamiltonian-connected graphs}

\author{ Petr Vr\'ana$^{1,4}$, \ Xingzhi Zhan$^{2,}$\thanks{Corresponding author.}, \  Leilei Zhang$^{2}$}

\date{\version}

\newcounter{mathitem}
\newenvironment{mathitem}
  {\begin{list}{{$(\roman{mathitem})$}}{
   \setcounter{mathitem}{0}
   \usecounter{mathitem}
   \setlength{\topsep}{0pt plus 2pt minus 0pt}
   \setlength{\parskip}{0pt plus 2pt minus 0pt}
   \setlength{\partopsep}{0pt plus 2pt minus 0pt}
   \setlength{\parsep}{0pt plus 2pt minus 0pt}
   \setlength{\leftmargin}{35pt}
   \setlength{\itemsep}{0pt plus 2pt minus 0pt}}}
  {\end{list}}

\newcounter{mylist}
  {\begin{list}{{$\roman{mylist}$}}{
   \setlength{\topsep}{7pt plus 2pt minus 0pt}
   \setlength{\parsep}{3pt plus 2pt minus 0pt}
   \setlength{\leftmargin}{12pt}
   \setlength{\labelwidth}{-6pt}
   }
   }
  {\end{list}}

\newcounter{prostredi}
\def\theprostredi{\arabic{prostredi}}
\def\vejde#1{\unskip
\nobreak\hfill\penalty50\hskip1em\hbox{}\nobreak\hfill
\hbox{#1}}
\newenvironment{theorem}{\par\bigskip\noindent%
\refstepcounter{prostredi}{\bf Theorem \theprostredi.}\quad\bgroup\sl }
{\egroup\par\bigskip\endtrivlist}%
{\egroup\vejde{\rule{2.5mm}{2.5mm}}\par\bigskip\endtrivlist}%
\newenvironment{proofbt}{\par
\noindent%
{\bf Proof}\bgroup}
{\egroup\vejde{\rule{2.5mm}{2.5mm}}\par\bigskip\endtrivlist}%
\newenvironment{corollary}{\par\bigskip\noindent%
\refstepcounter{prostredi}{\bf Corollary
\theprostredi.}\quad\bgroup\sl }
{\egroup\par\bigskip\endtrivlist}%
\newcounter{prostrclaim}
\def\theprostrclaim{\arabic{prostrclaim}}
\def\vejde#1{\unskip
\nobreak\hfill\penalty50\hskip1em\hbox{}\nobreak\hfill
\hbox{#1}}
{\egroup\vejde{$\Box$}\par\bigskip\endtrivlist}%
{\egroup\vejde{
$\Box$}\par\bigskip\endtrivlist}%

\newcounter{prostralph}
\def\theprostralph{\Alph{prostralph}}
\def\vejde#1{\unskip
\nobreak\hfill\penalty50\hskip1em\hbox{}\nobreak\hfill
\hbox{#1}}
\newenvironment{theoremAcite}[1]{\par\bigskip\noindent%
\refstepcounter{prostralph}{\bf Theorem
\theprostralph{} {#1}.}\quad\bgroup\sl }
{\egroup\par\bigskip\endtrivlist}%
\newenvironment{lemmaAcite}[1]{\par\bigskip\noindent%
\refstepcounter{prostralph}{\bf Lemma
\theprostralph{} {#1}.}\quad\bgroup\sl }
{\egroup\par\bigskip\endtrivlist}%
\newenvironment{conjectureAcite}[1]{\par\bigskip\noindent%
\refstepcounter{prostralph}{\bf Conjecture
\theprostralph{} {#1}.}\quad\bgroup\sl }
{\egroup\par\bigskip\endtrivlist}%

\newcommand{\noi}{\noindent}

\newcommand{\iso}{\simeq}

\newcommand{\cl}{{\rm cl}}

\newcommand{\co}{{\rm co}}

\newcommand{\cR}{{\cal R}}

\newcommand{\bp}{\beginpicture}
\newcommand{\ep}{\endpicture}

\newcommand{\Lp}{L^{-1}}

\newcommand{\ms}{\medskip}

\begin{document}
\maketitle

\footnotetext[1]{Department of Mathematics; European Centre of Excellence
  NTIS - New Technologies for the Information Society, University of West Bohemia,
  Univerzitn\'{\i}~8, 301 00 Pilsen, Czech Republic}
\footnotetext[2]{Department of Mathematics, East China Normal University, Shanghai 200241, China}
\footnotetext[3]{E-mail addresses: {vranap$@$kma.zcu.cz(P.Vr\'ana),} {zhan@math.ecnu.edu.cn(X.Zhan),} {mathdzhang@163.com(L.Zhang).}}
\footnotetext[4]{Research supported by project GA20-09525S of the Czech Science Foundation}

\begin{abstract}
\noi
We prove that every $3$-connected claw-free graph with domination number at most $3$ is hamiltonian-connected. The result is sharp and it is inspired by a conjecture posed by Zheng, Broersma, Wang and Zhang in 2020.

\ms

\noi

\noi
\end{abstract}

\section{Introduction}
\label{sec-intro}

In this paper, by a {\em graph} we always mean a simple finite undirected graph;
if we admit multiple edges, we always speak about a {\em multigraph}.
We follow the book \cite{W96} for terminology and notations.

The complete bipartite graph on $s$ and $t$ vertices is denoted by $K_{s,\, t}$.
The graph $K_{1,3}$ is called the {\it claw}. A graph is called {\it claw-free}
if it contains no induced subgraph isomorphic to the claw. A graph is called
{\it hamiltonian-connected} if between any two distinct vertices there is a hamiltonian
path. A subset $X$ of vertices in a graph $G$ is called a {\it dominating set} if
every vertex of $G$ is either contained in $X$ or adjacent to some vertex of $X$.
The {\em domination number} of $G$ is the size of a smallest dominating set of $G$.

In 1994, Ageev \cite{A94} proved the following sufficient condition for a claw-free
graph to be hamiltonian involving the domination number.

%
%
\begin{theoremAcite}{\cite{A94}}
\label{thmA-ham}
Every $2$-connected claw-free graph with domination number at most $2$ is hamiltonian.
\end{theoremAcite}

The main result of this note is inspired by the following conjecture  posed by Zheng, Broersma, Wang and Zhang \cite{ZBWZ}.
Note that a hamiltonian-connected graph is necessarily $3$-connected.

%
%
\begin{conjectureAcite}{\cite{ZBWZ}}
 Every $3$-connected claw-free graph with domination number at most $2$ is hamiltonian-connected.
\end{conjectureAcite}

We prove the following stronger theorem which is sharp.

%
%
\begin{theorem}
\label{thm1-3-dom}
Every 3-connected claw-free graph with domination number at most
3 is hamiltonian-connected.
\end{theorem}

The proof of Theorem \ref{thm1-3-dom} is postponed to Section \ref{sec-2}.
We will first need to recall some necessary known concepts and results.
%
%
We say that an edge is {\em pendant} if it contains a vertex of degree $1$,
and that a vertex is {\it simplicial} if its neighbors induce a complete graph.
The {\em line graph} of a multigraph $H$ is the graph $G=L(H)$ with $V(G)=E(H)$,
in which two vertices are adjacent if and only if the corresponding edges of $H$
have at least one vertex in common.

The following was proved in \cite{RV11-2} using a modification
of an approach from \cite{Z97}.

%
%
\begin{theoremAcite}{\cite{RV11-2}}
\label{thmA-vzor-jedn}
Let $G$ be a connected line graph of a multigraph. Then there is, up
to an isomorphism, a uniquely determined multigraph $H=\Lp(G)$ such that
a vertex $e\in V(G)$ is simplicial in $G$ if and only if the corresponding
edge $e\in E(H)$ is a pendant edge in $H$.
\end{theoremAcite}

An edge-cut $R\subset E(H)$ of a multigraph $H$ is {\em essential} if $H-R$ has
at least two nontrivial components, and $H$ is
{\em essentially $k$-edge-connected} if every essential edge-cut of $H$ is
of size at least $k$. It is a well-known fact that a line graph $G$ is
$k$-connected if and only if $\Lp(G)$ is essentially $k$-edge-connected.

A set of vertices $M\subset V(G)$ {\em dominates} an edge $e$ if $e$ has at
least one vertex in $M$. A~closed trail $T$ is a {\em dominating closed trail}
(abbreviated DCT) if $T$ dominates all edges of $G$ and an $(e,f)$-trail
(i.e, a trail with terminal edges $e,f$) is an {\em internally dominating
$(e,f)$-trail } (abbreviated $(e,f)$-IDT) if the set of its interior vertices
dominates all edges of $G$.

Harary and Nash-Williams \cite{HNW} proved a correspondence between a DCT
in $H$ and a hamiltonian cycle in $L(H)$ (the result was established in \cite{HNW}
for graphs, but it is easy to observe that the proof is true also for line graphs
of multigraphs).
A similar result showing that $G=L(H)$ is hamiltonian-connected if and only
if $H$ has an $(e_1,e_2)$-IDT for any pair of edges $e_1,e_2\in E(H)$, was
given in \cite{LLZ05}.

%
%
\begin{theoremAcite}{\cite{HNW,LLZ05}}
\label{thmA-DCT+IDT}
Let $H$ be a multigraph with $|E(H)|\geq 3$ and let $G=L(H)$.
\begin{mathitem}
\item {\bf\cite{HNW}} The graph $G$ is hamiltonian if and only if $H$ has a DCT.
\item {\bf\cite{LLZ05}}  For every $e_i\in E(H)$ and $a_i=L(e_i)$, $i=1,2$,
      $G$ has a hamiltonian $(a_1,a_2)$-path if and only if $H$ has an
      $(e_1,e_2)$-IDT.
\end{mathitem}
\end{theoremAcite}

Let $G$ be a 3-connected line graph and let $H = \Lp(G)$. The {\em  core}
of $H$ is the multigraph $\co(H)$ obtained from $H$ by removing all pendant
edges and suppressing all vertices of degree~2.

Shao {\cite{S05}} proved the following properties of the core of a multigraph.

%
%
\begin{lemmaAcite}{\cite{S05}}
\label{thmA-core}
Let $H$ be an essentially 3-edge-connected multigraph. Then
\begin{mathitem}
\item $\co(H)$ is uniquely determined,
\item $\co(H)$ is 3-edge-connected,
\item if $\co(H)$ has a spanning closed trail, then $H$ has a dominating closed trail.
\end{mathitem}
\end{lemmaAcite}

\begin{figure}[ht]
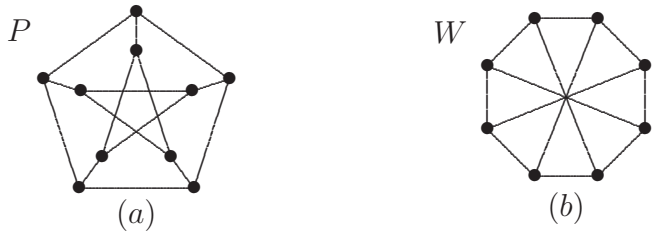

$$\beginpicture
\setcoordinatesystem units <0.7mm,1mm>
\setplotarea x from -70 to 70, y from -5 to 5
\put{\bp
\setcoordinatesystem units <.065mm,.065mm>
\setplotarea x from -200 to 200, y from -200 to 200
\put{$\bullet$} at  0 200
\put{$\bullet$} at 190.211303259 61.80339887499
\put{$\bullet$} at  117.5570504585  -161.803398875
\put{$\bullet$} at -117.5570504585  -161.803398875
\put{$\bullet$} at -190.211303259 61.80339887499
\plot 0 200 190.211303259 61.80339887499
        117.5570504585  -161.803398875 -117.5570504585
        -161.803398875 -190.211303259 61.80339887499 0 200 /
\plot  0 200 0 120 /
\plot 190.211303259 61.80339887499  114.1267819554 37.08203932499 /
\plot 117.5570504585  -161.803398875  70.5342302751  -97.082039325 /
\plot -117.5570504585  -161.803398875 -70.5342302751 -97.082039325 /
\plot -190.211303259 61.80339887499 -114.1267819554  37.08203932499 /
\setcoordinatesystem units <.039mm,.039mm>
\put{$\bullet$} at  0 200
\put{$\bullet$} at 190.211303259 61.80339887499
\put{$\bullet$} at  117.5570504585  -161.803398875
\put{$\bullet$} at -117.5570504585  -161.803398875
\put{$\bullet$} at -190.211303259 61.80339887499
\plot 0 200 117.5570504585  -161.803398875
    -190.211303259 61.80339887499 190.211303259 61.80339887499
    -117.5570504585 -161.803398875 0 200 /
\put{$P$}   at  -400    270
\put{$(a)$} at     0   -360
\ep} at -40  0
\put {\bp
\setcoordinatesystem units <0.35mm,0.35mm>
  \put{$\bullet$} at   -30   12
  \put{$\bullet$} at   -30  -12
  \put{$\bullet$} at    30   12
  \put{$\bullet$} at    30  -12
  \put{$\bullet$} at   -12   30
  \put{$\bullet$} at    12   30
  \put{$\bullet$} at   -12  -30
  \put{$\bullet$} at    12  -30
  \plot -30 -12 -30 12 -12 30 12 30 30 12 30 -12 12 -30
        -12 -30 -30 -12 /
  \plot -30  12  30 -12 /
  \plot -30 -12  30  12 /
  \plot -12  30  12 -30 /
  \plot -12 -30  12  30 /
%
%
%
%
 \put{$W$}   at  -44    25
 \put{$(b)$} at    0   -42
\ep} at 40 0
\endpicture$$
\vspace*{-6mm}
\caption{The Petersen graph $P$ and the Wagner graph $W$}
\label{fig-Petersen-Wagner}
\end{figure}

We denote by $P$ the Petersen graph and by $W$ the Wagner graph
(see Fig.~\ref{fig-Petersen-Wagner}).

%
%

Let $G$ be a multigraph, $R\subset G$ a spanning subgraph of $G$, and let $\cR$
be the set of components of $R$. Then $G/R$ is the multigraph with $V(G/R)=\cR$,
in which, for each edge in $E(G)$ between two components of $R$, there is an
edge in $E(G/R)$ joining the corresponding vertices of $G/R$.
The (multi-)graph $G/R$ is said to be a {\em contraction} of $G$. (Roughly,
in $G/R$, components of $R$ are contracted to single vertices while keeping
the adjacencies between them).

The contraction operation maps $V(G)$ onto $V(G/R)$ (where vertices of a
component of $R$ are mapped on a vertex of $G/R$). If $G/R\iso F$, then
this defines a function $\alpha: G\rightarrow F$ which is called
a {\em contraction of $G$ on $F$}.

The following theorem was proved in {\cite{LRVXY20+I}} (see also \cite{RV21}).

%
%
\begin{theoremAcite}{\cite{LRVXY20+I}}
\label{thmA-8vertices_plus_edge}
Let $H$ be a 3-edge-connected multigraph, $A\subset V(H)$, $|A|=8$, and let
$e\in E(H)$. Then either
\begin{mathitem}
\item $H$ contains a closed trail $T$ such that $A\subset V(T)$ and $e\in E(T)$,
      or
\item there is a contraction $\alpha:H\rightarrow P$ such that
      $\alpha(e)=xy\in E(P)$ and $\alpha(A)=V(P)\setminus\{x,y\}$.
\end{mathitem}
\end{theoremAcite}

In fact, we will need only the following easy corollary.

%
%
\begin{corollary}
\label{thmA-8vertices_plus_edge2}
Let $H$ be a 3-edge-connected multigraph, $A\subset V(H)$, $|A|\le 7$, and let
$e\in E(H)$. Then $H$ contains a closed trail $T$ such that $A\subset V(T)$ and $e\in E(T)$.
\end{corollary}

The concept of an M-closure $\cl^M(G)$ of a claw-free
graph $G$ was defined in \cite{RV11-2}. We do not need to know the exact construction 	of this closure. We will use only the following theorem proved in \cite{RV11-2}.

%
%
\begin{theoremAcite}{\cite{RV11-2}}
\label{thmA-strong_clos}
Let $G$ be a claw-free graph and let the graph $\cl^M(G)$ be its M-closure.
Then $\cl^M(G)$ has the following properties:
\begin{mathitem}
\item $V(G)=V(\cl^M(G))$ and $E(G)\subset E(\cl^M(G))$,
\item $\cl^M(G)$ is uniquely determined,
\item $G$ is hamiltonian-connected if and only if $\cl^M(G)$ is
         hamiltonian-connected,
\item[$(vi)$] $\cl^M(G)=L(H)$, where $H$ is a multigraph.
\end{mathitem}
\end{theoremAcite}

\section{Proof of Theorem~\ref{thm1-3-dom}}
\label{sec-2}

Now we are ready to give a proof of the main result.

\bigskip


\begin{proofbt} {\bf of Theorem~\ref{thm1-3-dom}.} \quad

Let $G$ be a counterexample to Theorem \ref{thm1-3-dom}.
By Theorem~\ref{thmA-strong_clos}, $\cl^M(G)$ is a non-hamiltonian-connected line
graph of a multigraph. Let $H=\Lp(\cl^M(G))$.
By Theorem \ref{thmA-DCT+IDT}, there are edges $e_1, e_2$ such that $H$ has
no $(e_1,e_2)$-IDT. Since $G$ is 3-connected, $H$ is essentially 3-edge-connected.

To reach a contradiction, we first convert the problem into the core of $H$, and
then we find a trail such that the corresponding trail in $H$ is an $(e_1,e_2)$-IDT.

For $i\in\left\lbrace 1,2 \right\rbrace,$ if $e_i,$ as an edge in $H,$ has both end-vertices of degree
at least 3, then $e_i\in E(\co(H))$, and we set $e^0_i=e_i.$
If $e_i$ is a pendant edge, then we denote by $e^0_i$ an arbitrary edge in $\co(H)$
containing the vertex of higher degree of $e_i$.
The last case is that one of the end-vertices of $e_i$, say, $v_2$, has degree 2.
In  this case we take as $e^0_i$ the new edge in $\co(H)$ resulting by suppressing
the vertex $v_2$.

Now, if $e^0_1=e^0_2$, we set $e_n=e^0_1=e^0_2$; otherwise we subdivide the edges
$e^0_1,e^0_2$ and join the two new vertices with an edge $e_n$.
In both cases, we denote by $H_n$ the resulting graph.
In the first case, $H_n=\co(H)$ which is $3$-edge-connected by Lemma \ref{thmA-core},
and it is easy to see that in the second case, $H_n$ is also $3$-edge-connected.

Let $\{w_1,w_2,w_3\}$ be a dominating set of $G$. By the definition of a line graph,
the three corresponding edges in $H$, denoted $f_1,f_2,f_3$, dominate all edges
of $H$. For $f_1,f_2,f_3\in E(H)$, we find edges $f^0_j,j=1,2,3$ in $\co(H)$ by
the same rules as $e^0_i$.
Since the three edges $f^0_j,j=1,2,3$ have at most 6 different vertices, denoted
$z_1,\ldots,z_6$, by Corollary~\ref{thmA-8vertices_plus_edge2}, there is a closed
trail  in $H_n$ containing the vertices $z_1,\ldots,z_6$ and the edge $e_n$.
It is straightforward to check that for every case of constructing $e^0_i$ and $e_n$,
we can find an $(e_1,e_2)$-IDT in $H$ (recall that the corresponding set of vertices
$\lbrace z_1,\ldots,z_6\rbrace$ dominates all edges of $H$).
\end{proofbt}

The result is sharp. A counterexample for domination number 4 is the line graph
of a graph obtained from the Wagner graph by adding at least one pendant edge
to each of its vertices.

\vskip 5mm
{\bf Acknowledgement.} This study was partly supported by project GA20-09525S of the Czech Science Foundation
(Petr Vr\'ana). This research  was also supported by the NSFC grants 11671148 and 11771148 and Science and Technology Commission of Shanghai Municipality (STCSM) grant 18dz2271000(Xingzhi Zhan, Leilei Zhang).

\end{document}